\documentclass[final,1p,times]{elsarticle}

\usepackage{amssymb}
\usepackage{amsmath}

\journal{Number Theory}
\begin{document}

\begin{frontmatter}


\title{On Taylor series expansion of $(1+ z)^{A}$ for $|z|>1$}

\author{Akhila Raman }

\address{University of California at Berkeley, CA-94720. Email: akhila.raman@berkeley.edu. Ph: 510-540-5544}



\begin{abstract}

It is well known that the Taylor series expansion of  $(1+ z)^{A}$ does not
converge for $|z|>1$ where A is a real number which is not equal to zero or a positive integer. 
A limited series expansion of this expression is obtained in this paper for 
$|z|>1$ as a product of convergent series.

\end{abstract}


\end{frontmatter}
   






\section{\label{sec:level1}Introduction\protect\\  \lowercase{} }

It is well known that the Taylor series expansion of  $(1+ z)^{A}$ is given by

\begin{equation} 
(1+ z)^{A}= \sum_{n=0}^{\infty} \binom{A}{n}  z^{n} \end{equation}

 for $|z|<1$ where  $\binom{A}{n}$ is the binomial choose function. It is well known that 
the series expansion does not converge for $|z|>1$ where A is a real number which is not equal to zero or a positive integer.

We could obtain a limited series expansion for $|z|>1$ by writing the above expression as follows

\begin{equation}    (1+ z)^{A}  =  (1+ \frac{z}{2} + \frac{z}{2}  )^{A} =  (1+ \frac{z}{2}  )^{A}  (1+ \frac{ \frac{z}{2}}{ (1+ \frac{z}{2}) }  )^{A}  = (1+ \frac{z}{2}  )^{A}  (1+ \frac{z}{ z+ 2 }  )^{A}    \end{equation}

The second term in the above equation has a convergent series representation, given that $|\frac{z}{ z+ 2 }| < 1$. If   $|\frac{z}{2}|>1$, we can write
\begin{equation}    (1+ \frac{z}{2})^{A}  =   (1+ \frac{z}{4}  )^{A}  (1+ \frac{z}{ z+ 4 }  )^{A}    \end{equation}

Repeating this procedure iteratively, if   $m_{0}$ is the minimum value for which $|\frac{z}{2^{m_{0}}}|<1$, we can write

\begin{equation}     (1+ z)^{A}  =   (1+ \frac{z}{2^{m_{0}}})^{A} \prod^{m_{0}}_{r=1}  (1+ \frac{z}{ z+   2^{r} }  )^{A}    \end{equation}

Each of the terms in the above product of terms has a convergent series representation. Given that we can write the convergent series expansion for each of the terms above as 
$(1+ \frac{z}{2^{m_{0}}})^{A}= \sum_{n=0}^{\infty} \binom{A}{n}  (\frac{z}{2^{m_{0}}})^{n}$
and $(1+ \frac{z}{ (z+   2^{r}) }  )^{A} =  \sum_{m=0}^{\infty} \binom{A}{m} (\frac{z}{ (z+   2^{r}) })^{m}$, where $\binom{A}{n}$ represents the Choose function[1],  we have the \textbf{series expansion for $(1+ z)^{A}$ expressed 
as a product of convergent series, which converges for } $|z|>1$  as follows:

\begin{equation}  (1+ z)^{A}   =  [ \sum_{n=0}^{\infty} \binom{A}{n}  (\frac{z}{2^{m_{0}}})^{n}     ] [ \prod^{m_{0}}_{r=1}  \sum_{m=0}^{\infty} \binom{A}{m}  [\frac{z}{ z+   2^{r}}]^{m} ] \end{equation}

\section{\label{sec:level1}Section 2\protect\\  \lowercase{} }

Let us take the case of $m_{0}=1$ for  $1<|z|<2$.  Expanding the  term $\frac{1}{(z+   2^{r})^{m}}$ in the above equation 5 as
Taylor series around a point $z=0$ , we have for $m>0$

\begin{equation}   \frac{1}{(z+   2^{r})^{m}} = \sum_{j=0}^{\infty} b(j,r,m) z^{j}    \end{equation}

where $b(0,r,m)= \frac{1}{(2^{r})^{m}} $ and $b(j,r,m)$ is given as follows for $j=1,2,3...$

\begin{equation}  b(j,r,m)= \binom{m+j-1}{j} \frac{(-1)^{j}}{(2^{r})^{m+j}} ; \end{equation}

For $m=0$, $\frac{1}{(z+   2^{r})^{m}}=1$. Now we can write the the series expansion of $(1+ z)^{A}$ which converges for $1<|z|<2$, as a product of terms expanded in Taylor series  as follows:

\begin{equation}  (1+ z)^{A}  =  [ \sum_{n=0}^{\infty} \binom{A}{n}  (\frac{z}{2})^{n} ]  [ \sum_{m=0}^{\infty} \binom{A}{m}   z^{m}   \sum_{j=0}^{\infty} b(j,1,m) z^{j} ] \end{equation}

For the case of $m_{0}>1$ for $|z|>2$, we can write as follows:

\begin{equation}  (1+ z)^{A}  =  [ \sum_{n=0}^{\infty} \binom{A}{n}  (\frac{z}{2^{m_{0}}})^{n} ]  [ \sum_{m=0}^{\infty} \binom{A}{m}   z^{m}   \sum_{j=0}^{\infty} b(j,m_{0},m) z^{j} ]  [ \prod^{m_{0}-1}_{r=1}  \sum_{m=0}^{\infty} \binom{A}{m} z^{m}  [\frac{1}{ z+   2^{r}}]^{m} ] \end{equation}

The last term in the above equation $ (\frac{1}{ z+   2^{r}})^{m} $ 
can be expressed as follows:

\begin{equation} ( \frac{1}{ z+   2^{r}})^{m}  = (z+   2^{r})^{-m} = 2^{-r*m}  (1+ \frac{z}{2^{r}})^{-m}  \end{equation}

The term$(1+ \frac{z}{2^{r}})^{-m}$ can be  recursively expanded using Eq.9 by substituting $z \to \frac{z}{2^{r}}$ and $A \to -m$ 
and $m_{0} \to m_{0}-r$ to obtain  the series expansion of 
$(1+ z)^{A}$ which converges for $|z|>1$ as a product of terms expanded in Taylor series.

\section{\label{sec:level1}Section 3\protect\\  \lowercase{} }

Let us consider the following binomial expression  \\

\begin{equation}  (x + y )^{A}   \end{equation}

where $A$ is a real number which is not equal to zero or a positive integer and $z=\frac{x}{y}$ and $|z|>1$. Writing $(x + y )^{A} = (1 + \frac{x}{y} )^{A} y^{A} = (1 + z )^{A} y^{A} $, we can write the series expansion of this  expression using results obtained in equations 5 and 9 as follows:\\

\begin{equation}  (x + y)^{A}   =  y^{A}  \sum_{n=0}^{\infty} \binom{A}{n}  (\frac{z}{2^{m_{0}}})^{n}      \prod^{m_{0}}_{r=1}  \sum_{m=0}^{\infty} \binom{A}{m}  [\frac{z}{ z+   2^{r}}]^{m}  \end{equation}

\begin{equation}  (x + y)^{A}  = y^{A} [ \sum_{n=0}^{\infty} \binom{A}{n}  (\frac{z}{2^{m_{0}}})^{n} ]  [ \sum_{m=0}^{\infty} \binom{A}{m}   z^{m}   \sum_{j=0}^{\infty} b(j,m_{0},m) z^{j} ]  [ \prod^{m_{0}-1}_{r=1}  \sum_{m=0}^{\infty} \binom{A}{m} z^{m}  [\frac{1}{ z+   2^{r}}]^{m} ]  \end{equation}

\section{\label{sec:level1} Conclusion \protect\\  \lowercase{} }

It has been shown that the Taylor series expansion of  $(1+ z)^{A}$ can be expanded
as a product of convergent series, for $|z|>1$ 
where A is a real number which is not equal to zero or a positive integer. 

\section{\label{sec:level1} Acknowledgements \protect\\  \lowercase{} }

The author would like to thank Michael Schlosser(Universitat Wien) for his constructive comments on the paper.

\section{\label{sec:level1}References\protect\\  \lowercase{} }

[1]  Abramowitz, M. and Stegun, I. A. (Eds.).  Handbook of Mathematical Functions with Formulas, Graphs, and Mathematical Tables, 9th printing. New York: Dover,  1972.


\bibliographystyle{elsarticle-num}
\bibliography{<your-bib-database>}







\end{document}